\documentclass[a4paper, oneside]{amsart}
\usepackage{graphicx} 
\usepackage{amsmath}
\usepackage{amssymb}
\usepackage{amsfonts}
\usepackage{bbm}
\usepackage{mathtools}
\usepackage{cancel}
\usepackage[dvipsnames]{xcolor}

\usepackage[pagebackref=true,colorlinks=true, linkcolor=blue, citecolor=blue]{hyperref}
\renewcommand*\backref[1]{\ifx#1\relax \else (Cited on: p. #1) \fi}
\usepackage[nameinlink]{cleveref}

\newtheorem{theorem}{Theorem}[section]
\newtheorem{definition}[theorem]{Definition}

\newtheorem{lemma}[theorem]{Lemma}

\newtheorem{remark}[theorem]{Remark}

\newtheorem{notation}[theorem]{Notation}

\newtheorem*{proposition*}{Proposition}
\newtheorem*{theorem*}{Theorem}
\newtheorem*{criterion*}{Criterion}
\newtheorem{remark*}{Remark}

\DeclareMathOperator{\cov}{cov}

\newcommand{\connected}[3]{#1 \xleftrightarrow[]{#3} #2}
\newcommand{\notconnected}[3]{\cancel{#1 \xleftrightarrow[]{#3} #2}}
\newcommand{\g}{\mathfrak{g}}

\mathtoolsset{showonlyrefs=true}

\title[A correlation bound for the 1D Ising model]{A correlation bound for the one-dimensional heterogeneous Ising model}
\author{Edward Athaide}
\address{Department of Physics and Astronomy, University of California, Los Angeles.
}
\email{athaideej@physics.ucla.edu}
\author{Maciej G\l uchowski}
\address{Department of Mathematics, University of California, Los Angeles.
}
\email{mgluchowski@math.ucla.edu}
\author{Jonas K\"{o}ppl}
\address{Weierstrass Institute for Applied Analysis and Stochastics, Berlin.}
\email{koeppl@wias-berlin.de}
\author{Georg Menz}
\address{Department of Mathematics, University of California, Los Angeles.}
\email{gmenz@math.ucla.edu}  

\date{\today}

\numberwithin{equation}{section}

\begin{document}

\begin{abstract}
    We derive a new upper bound for the correlations in a heterogeneous one-dimensional Ising model with free boundary conditions. The new upper bound quantifies the simultaneous decay of correlations due to weakness of nearest-neighbor coupling constants and the effect of external fields. The proof constitutes an application of random currents to a non-ferromagnetic Ising model.
\end{abstract}

\maketitle

\section{Introduction}

\subsection{The model} 

We consider the one-dimensional heterogeneous Ising model on a finite subset $\{0, \ldots, N\} \subset \mathbb{Z}$, with nearest-neighbor edges. The model takes as parameters a set of two-site coupling coefficients $\mathbf{J} = \{J_x \in \mathbb{R}_{\geq 0}\}_{x=0}^{N-1}$ as well as site-specific external field strengths $\mathbf{h} = \{h_x \in \mathbb{R}\}_{x=0}^N$. Note that we require, without loss of generality, that all $J_x > 0$. The Hamiltonian for the model can be written as
\begin{equation}
    H_{\mathbf{J},\mathbf{h}} (\boldsymbol{\sigma}) := - \sum_{x=0}^{N-1} J_x \sigma_x \sigma_{x+1} - \sum_{x=0}^N h_x \sigma_x,
\end{equation}
where $\boldsymbol{\sigma} = \{ \sigma_x \}_{x=0}^N \in \{\pm 1\}^{N+1}$ is a spin configuration on the vertices of the lattice. We take $\boldsymbol{\sigma}$ to be a random variable distributed according to the Gibbs measure 
\begin{equation}
    \mu_{\mathbf{J},\mathbf{h}}(\boldsymbol{\sigma}) := \frac{1}{Z_{\mathbf{J},\mathbf{h}}} \exp(-H_{\mathbf{J},\mathbf{h}}(\boldsymbol{\sigma})), \label{equ:gibbs-measure}
\end{equation}
where we have used for normalization the partition function
\begin{equation}
    Z_{\mathbf{J},\mathbf{h}} := \sum_{\boldsymbol{\sigma} \in \{\pm 1\}^{N+1}} \exp(-H_{\mathbf{J},\mathbf{h}}(\boldsymbol{\sigma})).
\end{equation}
We have omitted a temperature parameter, which can be considered to be absorbed by the parameters $\mathbf{J}$ and $\mathbf{h}$. Expectations with respect to the Gibbs measure will be denoted by $\left< \cdot \right>_{\mathbf{J},\mathbf{h}}$, and covariances will be denoted by $\cov_{\mathbf{J},\mathbf{h}}(\cdot, \cdot)$.

If a particular subinterval $\{i, \ldots, j\}$ of the lattice is of interest, then integrating out spin variables situated on $\{0, \ldots, i-1\} \cup \{j+1, \ldots, N\}$ yields a measure of a similar form to the original, with adjusted effective fields at points $i$ and $j$ to account for the influence of external spins. The following lemma makes this statement precise, and its proof is given in Appendix \ref{appendix:truncated-model}.
\begin{lemma}
    For given points $i < j$ on the lattice $\{0, \ldots, N\}$ with specified parameters $\mathbf{J}$ and $\mathbf{h}$, there exist effective fields $h'_i, h'_j \in \mathbb{R}$ at $i$ and $j$, respectively, such that the Gibbs measure on $(\sigma_i, \ldots, \sigma_j)$ defined as in Equation \eqref{equ:gibbs-measure} from the effective Hamiltonian
    \begin{equation}
        H'_{\mathbf{J},\mathbf{h'}}((\sigma_i, \ldots, \sigma_j)) := - \sum_{x=i}^{j-1} J_x \sigma_x \sigma_{x+1} - \left(  h'_i \sigma_i + \sum_{x=i+1}^{j-1} h_x \sigma_x + h'_j \sigma_j \right)
    \end{equation}
    is equal to the marginal measure of $\mu$ on $(\sigma_i, \ldots, \sigma_j)$. Moreover, if $i \neq 0$ then $|h'_i - h_i| \leq J_{i-1}$, and if $j \neq N$ then $|h'_j - h_j| \leq J_j$.
    \label{lem:effective-field}
\end{lemma}

\subsection{Decay of correlations of the 1D Ising model}

Of primary interest to us will be the quantity~$\cov_{\mathbf{J},\mathbf{h}}(\sigma_i, \sigma_j)$. For the homogeneous Ising model, one can apply the transfer matrix method introduced by Kramers and Wannier in \cite{KW41} to compute the covariance directly from the eigenvalues of a $2 \times 2$ transfer matrix product. If the coefficients are taken to be heterogeneous or random, then the transfer matrix method can still be used, with the matrix whose eigenvalues determine the covariance being a product of non-identical and possibly random matrices. The setting in which coefficients are i.i.d.~random has been particularly well-studied, starting with the work of Derrida and Hilhorst in \cite{DH81}, and many questions regarding the distribution of covariances for finite lattices in this setting are addressed in \cite{CNPV90,CPV12}. In the separate setting of a one-dimensional torus with no external field, some correlation inequalities are derived in \cite{CU08}, although these inequalities are concerned with the signs of covariances rather than estimating their magnitude.

Also of interest is the rate of exponential decay of the covariance for long distances; that is, the quantity
\begin{equation}
    \frac{1}{\xi} := \limsup_{\{i, \ldots, j\} \uparrow \mathbb{Z}} \left[ - \log \frac{\cov_{\mathbf{J},\mathbf{h}}(\sigma_i, \sigma_j)}{|j-i|} \right],
\end{equation}
where the limit is taken in a suitable sense. If the model parameters are selected i.i.d.~from a well-behaved distribution, then the above rate exists as a limit as the random transfer matrix products converge almost surely to deterministic growth rates \cite{Fur63,Ose65,Kin73,Lig85}. In general, the decay rate cannot be expressed in closed form. However, efficient algorithms for computing its value numerically have existed for several decades \cite{PS92,PPS95}. Moreover, the analytical dependence of the decay rate on properties of the distribution of model parameters has been studied in detail \cite{DH83,Hav19,GG21,CGGH25}. 

Nevertheless, the connection between this decay rate and the original parameters of the model is typically opaque based on existing literature. One simple bound in terms of the coupling constants results from the Ding-Song-Sun inequality \cite{DSS23}, which bounds the covariance by the corresponding covariance in a lattice under no external field, in which case the Gibbs measure decomposes into a product measure with $\cov_{\mathbf{J},\mathbf{0}}(\sigma_i,\sigma_j) = \prod_{x=i}^{j-1} \tanh(J_x)$. An exponential decay rate in the presence of constant external fields was found in \cite{LO68}, and an argument for a similar result presented in \cite{Iof09} served as an inspiration for the arguments presented here. Our results integrate these two conclusions into a single estimate.

\subsection{Summary of Results}

Our main result is the following upper bound on covariances. 
\begin{theorem}
    For a set of nonnegative coupling constants $\mathbf{J}$, any external field $\mathbf{h}$, and any two points $i<j$ on the lattice $\{0, \ldots, N\}$, we have
    \begin{equation}
        \cov_{\mathbf{J},\mathbf{h}} (\sigma_i, \sigma_j) \leq \left( \prod_{x=i}^{j-1} \frac{4 \tanh(J_x)}{(1+\tanh(J_x))^2} \right) \frac{4 e^{-2 \left|h'_i + \sum_{x=i+1}^{j-1} h_x + h'_j\right|}}{\left(1+e^{-2 \left( |h'_i| + \sum_{x=i+1}^{j-1} |h_x| + |h'_j| \right)} \right)^2},
    \end{equation}
    where $h'_i$ and $h'_j$ refer to the effective fields at points $i$ and $j$, respectively, obtained by integrating out all spin sites not between $i$ and $j$.
    \label{thm:cov-upper-bd-signed}
\end{theorem}
Theorem~\ref{thm:cov-upper-bd-signed} is derived from Theorem~\ref{thm:cov-upper-bd-fm}, an upper bound on covariances applying only to nonnegative external fields, and Lemma~\ref{lem:cov-signed-and-absolute}, an inequality relating Ising covariances due to signed interactions to their counterparts found when taking the absolute values of all interactions. We give both auxiliary statements below.

\begin{theorem}
    For a set of nonnegative coupling constants $\mathbf{J}$ and nonnegative external field $\mathbf{h}$, and any two points $i<j$ on the lattice $\{0, \ldots, N\}$,
    \begin{equation}
        \cov_{\mathbf{J},\mathbf{h}} (\sigma_i, \sigma_j) \leq \left( \prod_{x=i}^{j-1} \frac{4 \tanh(J_x)}{(1+\tanh(J_x))^2} \right) \frac{1}{\cosh^2 \left( {h'_i + \sum_{x=i+1}^{j-1} h_x + h'_j} \right)}.
    \end{equation}
    \label{thm:cov-upper-bd-fm}
\end{theorem}

\begin{lemma}
    For an Ising model on a lattice $\Lambda$ with $\mathbf{J}$ and $\mathbf{h}$ denoting respectively the two-site coupling coefficients and site-dependent external fields, both possibly signed, we have
    \begin{equation}
        |\cov_{\mathbf{J},\mathbf{h}}(\sigma_i, \sigma_j)| \leq \cov_{\mathbf{|J|},|\mathbf{h}|}(\sigma_i, \sigma_j) \left( \frac{Z_{\mathbf{|J|}, |\mathbf{h}|}}{Z_{\mathbf{J},\mathbf{h}}} \right)^2. \label{equ:cov-signed-and-absolute}
    \end{equation}
    for any two sites $i$ and $j$.
    \label{lem:cov-signed-and-absolute}
\end{lemma}

Both Theorem~\ref{thm:cov-upper-bd-fm} and Lemma~\ref{lem:cov-signed-and-absolute} take advantage of the random current expansion. The use of random currents to study covariances as in Theorem~\ref{thm:cov-upper-bd-fm} is standard, as in \cite{ABF87,ADS15,DT16}. The proof of Lemma~\ref{lem:cov-signed-and-absolute} relies on a recent extension of random currents to handle both positive and negative interactions introduced in \cite{Aiz25}. Such an application is atypical, as most past applications of random currents have handled only nonnegative interactions. An introduction to the random current expansion and proofs of Theorem~\ref{thm:cov-upper-bd-signed}, Theorem~\ref{thm:cov-upper-bd-fm}, and Lemma~\ref{lem:cov-signed-and-absolute} will be given in Section \ref{sec:main-results}.


\subsection{A note on prior versions of the article}

A previous version of this article involved a different handling of the extension of the random current representation to non-ferromagnetic interactions from what is presented here. While that version was in review, a similar, but in our view more powerful, treatment of non-ferromagnetic interactions was introduced in \cite{Aiz25}. In light of this development and reviewer feedback, we have restructured the article to adopt the conventions of \cite{Aiz25} and to present the proof of the main theorem as straightforwardly as possible. We recommend that the interested reader refer to \cite[Section~7]{Aiz25} for details on the use of random currents for non-ferromagnetic interactions.

\section{Proofs of the Main Results} \label{sec:main-results}

To prove Theorem \ref{thm:cov-upper-bd-fm} and Lemma \ref{lem:cov-signed-and-absolute}, we will use the random current expansion of the Ising model as introduced in \cite{Aiz82,Aiz25}. To begin this section, we will state some important results on random currents. Following this, we will outline and prove Theorem \ref{thm:cov-upper-bd-fm}. The section will conclude with proofs of Lemma \ref{lem:cov-signed-and-absolute} and Theorem \ref{thm:cov-upper-bd-signed}.

\subsection{The random current expansion for signed interactions} \label{ssc:random-currents}
The random current expansion is a tool for representing correlations between spins at different sites of an Ising model. We present it here allowing for general signed interactions. Taking a graph $\Lambda$ with edges $E \subset \Lambda \times \Lambda$, we can consider an Ising model with coupling constants $\{J_e \in \mathbb{R} \}_{e \in E}$ and no external field. The following definitions will provide the foundation for the discussion that follows:
\begin{definition}
    A current is a nonnegative integer-valued function $\mathbf{n}: E \to \mathbb{N}_0$. Given the coupling constants $\{J_e  \in \mathbb{R}\}_{e \in E}$, a random current refers to a current whose value is sampled from the measure
    \begin{equation}
        \mathbb{P}(\{\mathbf{n}\}) := \prod_{e \in E} \exp(-|J_e|) \frac{ |J_e|^{\mathbf{n}(e)}}{\mathbf{n}(e)!}. \label{equ:rc-measure}
    \end{equation}
    In other words, the random current measure $\mathbb{P}$ is a product of Poisson measures on each of the edges, with $|J_e|$ giving the rate of arrivals of $\mathbf{n}(e)$ on any edge $e \in E$.
\end{definition}

\begin{definition}
    For a (fixed) current $\mathbf{n}$, we define the boundary $\partial \mathbf{n}$ to be
    \begin{equation*}
        \partial \mathbf{n} := \left\{ x \in \Lambda: \sum_{e \in E: x \, \mathrm{touches} \, e} \mathbf{n}(e) \; \mathrm{is} \, \mathrm{odd} \right\}.
    \end{equation*}
\end{definition}

Because the boundary of a current will depend on the evenness or oddness of values that it takes on the edges, we recall the following facts about Poisson random variables.
\begin{remark}
    Let $X$ be a Poisson random variable with rate $\lambda \geq 0$. Then,
    \begin{align}
        \mathbb{P}_X(X = 0) &= \exp(-\lambda); \\
        \mathbb{P}_X(X \, \mathrm{even}) &= \exp(-\lambda) \cosh(\lambda); \\
        \mathbb{P}_X(X \, \mathrm{odd}) &= \exp(-\lambda) \sinh(\lambda).
    \end{align}
    \label{rmk:poisson-rv}
\end{remark}

It can be shown using arguments in \cite[Sections~4,7]{Aiz25} that for any set of vertices $A \subseteq \Lambda$, thermal moments of the spin variable take the form
\begin{equation}
    \Tilde{Z}_{\mathbf{J},\mathbf{0}} \left< \prod_{x \in A} \sigma_x \right>_{\mathbf{J},\mathbf{0}}  = \sum_{\mathbf{n}: \partial \mathbf{n} = A} \mathbb{P}(\{\mathbf{n}\}) (-1)^{\#^-_{\mathbf{J}}(\mathbf{n})} = \mathbb{E} \left[ (-1)^{\#^-_{\mathbf{J}}(\mathbf{n})} \mathbbm{1}_{\{\partial \mathbf{n} = A\}} \right], \label{equ:rc-sign-expectation}
\end{equation}
where we have defined
\begin{equation}
    \Tilde{Z}_{\mathbf{J},\mathbf{0}}  := \frac{Z_{\mathbf{J},\mathbf{0}} }{2^{|\Lambda|} \exp \left( \sum_{e \in E} |J_e| \right)},
\end{equation}
and $\#^-_{\mathbf{J}}(\mathbf{n})$ represents the total number of arrivals of $\mathbf{n}$ on edges where $\mathbf{J}$ is negative.

External fields are handled with Griffiths' ghost spin trick, first introduced in \cite{Gri67}, in which we append another site $\g$ to the lattice and fix $\sigma_{\g} = +1$ to break the global spin-flip symmetry. Along with $\g$, we add an additional edge to the graph connecting $\g$ to each vertex in $\Lambda$. Then, treating the external field $h_x$ at a site $x \in \Lambda$ as the coupling constant of an interaction between $x$ and $\g$, we can carry out the same expansion as above to obtain the result
\begin{equation}
    \Tilde{Z}_{\mathbf{J},\mathbf{h}} \left< \prod_{x \in A} \sigma_x \right>_{\mathbf{J},\mathbf{h}}  = \mathbb{E} \left[ (-1)^{\#^-_{\mathbf{J},\mathbf{h}}(\mathbf{n})} \mathbbm{1}_{\{\partial \mathbf{n} = A\} \cup \{\partial \mathbf{n} = A \cup \{\g\}\}} \right],
    \label{equ:rc-moments}
\end{equation}
where $\Tilde{Z}_{\mathbf{J},\mathbf{h}} $ and $\#_{\mathbf{J},\mathbf{h}}(\mathbf{n})$ are defined in the same way as $\Tilde{Z}_{\mathbf{J},\mathbf{0}} $ and $\#_{\mathbf{J}}(\mathbf{n})$, treating connections from sites in $\Lambda$ to $\g$ as additional edges.

A combinatorial argument from \cite[Lemma~1]{GHS70} when applied to the random current measure yields the so-called ``switching lemma,'' which was first stated in the context of random currents in \cite[Lemma~3.2]{Aiz82}. A particular consequence of that result that is of interest to us is the following:
\begin{lemma}[Consequence of the switching lemma]
    The covariance between two sites $i, j \in \Lambda$ can be written in the form
    \begin{equation}
        (\Tilde{Z}_{\mathbf{J},\mathbf{h}})^2 \cov_{\mathbf{J},\mathbf{h}} (\sigma_i, \sigma_j) = \mathbb{E}\left[ (-1)^{\#^-_{\mathbf{J},\mathbf{h}}(\mathbf{n}_1+\mathbf{n}_2)}\mathbbm{1}_{\{\partial \mathbf{n}_1 = \varnothing, \partial \mathbf{n}_2 = \{i, j\}, \notconnected{i}{\g}{\mathbf{n}_1 + \mathbf{n}_2}\}} \right], \label{equ:switching-lemma-alt}
    \end{equation}
    where the expectation is taken over two i.i.d.~random currents $\mathbf{n}_1, \mathbf{n}_2$, and the notation $\notconnected{i}{\g}{\mathbf{n}_1 + \mathbf{n}_2}$ denotes the event that there is no path from $i$ to $\g$ along edges on which either $\mathbf{n}_1$ or $\mathbf{n}_2$ takes a nonzero value.

    In particular, when $\mathbf{J}$ and $\mathbf{h}$ are all nonnegative, we have
    \begin{equation}
        \cov_{\mathbf{J},\mathbf{h}}(\sigma_i, \sigma_j) = \frac{\mathbb{P}(\{\partial \mathbf{n}_1 = \varnothing, \partial \mathbf{n}_2 = \{i, j\}, \notconnected{i}{\g}{\mathbf{n}_1 + \mathbf{n}_2}\})}{\left( \mathbb{P}(\{\partial \mathbf{n} = \varnothing\})\right)^2}. \label{equ:switching-lemma}
    \end{equation}
    \label{lem:switching-lemma}
\end{lemma}

\begin{proof}
    It follows from Equation~\eqref{equ:rc-moments} that
    \begin{equation}
        (\Tilde{Z}_{\mathbf{J},\mathbf{h}})^2 \left< \sigma_i \sigma_j \right> = \mathbb{E}\left[ (-1)^{\#^-_{\mathbf{J},\mathbf{h}}(\mathbf{n}_1+\mathbf{n}_2)}\mathbbm{1}_{\{\partial \mathbf{n}_1 = \varnothing, \partial \mathbf{n}_2 = \{i, j\}\}} \right]
    \end{equation}
    and
    \begin{equation}
        (\Tilde{Z}_{\mathbf{J},\mathbf{h}})^2 \left< \sigma_i \right> \left< \sigma_j \right> = \mathbb{E}\left[ (-1)^{\#^-_{\mathbf{J},\mathbf{h}}(\mathbf{n}_1+\mathbf{n}_2)}\mathbbm{1}_{\{\partial \mathbf{n}_1 = \{i,\g\}, \partial \mathbf{n}_2 = \{j, \g\}\}} \right].
    \end{equation}
    As the difference between these expressions will yield the left-hand side of Equation~\eqref{equ:switching-lemma-alt}, it suffices to show that
    \begin{equation}
        \begin{split}
            \mathbb{E}&\left[ (-1)^{\#^-_{\mathbf{J},\mathbf{h}}(\mathbf{n}_1+\mathbf{n}_2)}\mathbbm{1}_{\{\partial \mathbf{n}_1 = \{i,\g\}, \partial \mathbf{n}_2 = \{j, \g\}\}} \right] \\&= \mathbb{E}\left[ (-1)^{\#^-_{\mathbf{J},\mathbf{h}}(\mathbf{n}_1+\mathbf{n}_2)}\mathbbm{1}_{\{\partial \mathbf{n}_1 = \varnothing, \partial \mathbf{n}_2 = \{i, j\}, \connected{i}{\g}{\mathbf{n}_1 + \mathbf{n}_2}\}} \right].
        \end{split}
    \end{equation}
    For a proof of this fact, we refer the reader to \cite[Lemma~4.2]{Aiz25}.
\end{proof}

\subsection{Proof of Theorem~\ref{thm:cov-upper-bd-fm}} In this section, for simplicity we set $i=0$ and $j=N$, which by Lemma~\ref{lem:effective-field} extends to general $i, j$ making the proper substitutions $h_i \to h'_i$ and $h_j \to h'_j$. For notational convenience, we make some additional definitions. First, we decompose currents on the extended one-dimensional lattice with a ghost site into currents on the original edges of the lattice and currents on the edges appended to connect to the ghost site.
\begin{notation}
    In this section, $\mathcal{C}$ will denote the space of all currents on the extended lattice, $\mathcal{C}^{\Lambda}$ will denote the space of currents supported on nearest-neighbor edges within $\Lambda := \{0, \ldots, N\}$, and $\mathcal{C}^{\g}$ will denote the space of currents supported on edges connecting to the ghost site $\g$ introduced in Subsection \ref{ssc:random-currents}. Any $\mathbf{n} \in \mathcal{C}$ on the extended lattice with a ghost spin can be decomposed uniquely in the form $\mathbf{n} = \mathbf{n}^{\Lambda} \oplus \mathbf{n}^{\g}$, where $\mathbf{n}^{\Lambda} \in \mathcal{C}^{\Lambda}$ and $\mathbf{n}^{\g} \in \mathcal{C}^{\g}$. Moreover, $\mathbf{n}^{\g}$ and $\mathbf{n}^{\Lambda}$ are independent when $\mathbf{n}$ is sampled from $\mathbb{P}$. 
    \label{not:rc-decomposition}
\end{notation}
Additionally, we will use the following conventions:
\begin{notation}
    For a given $\mathbf{n}^{\g}$, we let $\mathbf{n}^{\g}(x)$ for $x \in \{0, \ldots, N\}$ denote the value taken by $\mathbf{n}^{\g}$ on the edge connecting $x$ to $\g$. Similarly, for a given $\mathbf{n}^{\Lambda}$, we let $\mathbf{n}^{\Lambda}(x)$ for $x \in \{0, \ldots, N-1\}$ denote the value taken by $\mathbf{n}^{\Lambda}$ on the edge $(x, x+1)$. In general, edges will be indexed by their left vertex.
\end{notation}

To prove Theorem~\ref{thm:cov-upper-bd-fm}, we first show that Equation \eqref{equ:switching-lemma} in the one-dimensional context, with currents decomposed as in Notation \ref{not:rc-decomposition}, yields the following formula:
\begin{lemma}
    The spin variables at the endpoints of $\{0, \ldots, N\}$ obey the covariance bound
    \begin{equation}
        \cov(\sigma_0, \sigma_N) = \prod_{x=0}^{N-1} \tanh(J_x) \left[ \frac{\mathbb{P}(\{\mathbf{n}^{\Lambda} \, \mathrm{even}\}) \mathbb{P}(\{\mathbf{n}^{\g} \equiv 0\})}{\mathbb{P}\left( \{\partial \mathbf{n}^{\Lambda} = \partial \mathbf{n}^{\g}\} \right)} \right]^2. \label{equ:cov-identity}
    \end{equation}
    \label{lem:cov-identity}
\end{lemma}
The main simplification due to the constrained nature of the one-dimensional setting occurs in the numerator of Equation \eqref{equ:switching-lemma}. 

\begin{lemma}
    The event $\{\partial \mathbf{n}_1 = \varnothing, \partial \mathbf{n}_2 = \{0, N\}, \notconnected{0}{\g}{\mathbf{n}_1 + \mathbf{n}_2}\}$ is equivalent to the following three conditions all being met:
    \begin{enumerate}
        \item[(i)] $\mathbf{n}_1^{\Lambda}(x)$ is even for all $x \in \{0, \ldots, N-1\}$;
        \item[(ii)] $\mathbf{n}_2^{\Lambda}(x)$ is odd for all $x \in \{0, \ldots, N-1\}$; and
        \item[(iii)] $\mathbf{n}_1^{\g}$ and $\mathbf{n}_2^{\g}$ are both identically zero.
    \end{enumerate}
    \label{lem:cov-numerator-1d}
\end{lemma}

\begin{proof}
    First, suppose that $(\mathbf{n}_1,\mathbf{n}_2) \in \mathcal{C}^2$ belongs to the event in the numerator of Equation \eqref{equ:switching-lemma}. The condition that $\notconnected{0}{\g}{\mathbf{n}_1+\mathbf{n}_2}$ implies that $\mathbf{n}_1^{\g}(0) = \mathbf{n}_2^{\g}(0) = 0$, and $\mathbf{n}_1^{\g}(x) = \mathbf{n}_2^{\g}(x) = 0$ for any $x$ that is connected to $0$ via $\mathbf{n}_1 + \mathbf{n}_2$. 
    
    Using these observations and considering the requirements that $0 \notin \partial \mathbf{n}_1$ and that $0 \in \partial \mathbf{n}_2$, we deduce that $\mathbf{n}_1^{\Lambda}(0)$ is even and that $\mathbf{n}_2^{\Lambda}(0)$ is odd. Thus, vertex $1$ is connected to $0$ at least by $\mathbf{n}_2$, and so $\mathbf{n}_1^{\g}(1) = \mathbf{n}_2^{\g}(1) = 0$. Proceeding inductively, if vertex $x > 0$ is connected to $0$ with $\mathbf{n}_1^{\Lambda}(x-1)$ being even and $\mathbf{n}_2^{\Lambda}(x-1)$ being odd, then $\mathbf{n}_1^{\g}(x) = \mathbf{n}_2^{\g}(x) = 0$, and if $x < N$ then the fact that $x \notin \partial \mathbf{n}_1 \cup \partial \mathbf{n}_2$ implies that $\mathbf{n}_1^{\Lambda}(x)$ is even while $\mathbf{n}_2^{\Lambda}(x)$ is odd and connects $x+1$ to all vertices to its left. This recursion continues until it reaches $j$, verifying that the conditions listed in Lemma \ref{lem:cov-numerator-1d} are all met.

    Conversely, suppose that all three conditions (i)--(iii) are satisfied. By examining the outgoing currents at any vertex, it is easy to verify that $\partial \mathbf{n}_1 = \varnothing$ and that $\partial \mathbf{n}_2 = \{0,N\}$. Moreover, given that $\mathbf{n}_1^{\g}$ and $\mathbf{n}_2^{\g}$ are identically zero, no points $0, \ldots, N$ connect to $\g$, and thus $(\mathbf{n}_1, \mathbf{n}_2)$ belongs to the event in the numerator of \eqref{equ:switching-lemma}.
\end{proof}

The proof of Lemma~\ref{lem:cov-identity} follows:
\begin{proof}[Proof of Lemma~\ref{lem:cov-identity}]
    As a result of Lemma \ref{lem:cov-numerator-1d}, we can write
    \begin{equation}
        \begin{split}
            \mathbb{P}&(\{\partial \mathbf{n}_1 = \varnothing, \partial \mathbf{n}_2 = \{0, N\}, \notconnected{0}{\g}{\mathbf{n}_1 + \mathbf{n}_2}\}) \\ &= \mathbb{P}(\{\mathbf{n}_1^{\Lambda} \, \text{even}\}) \mathbb{P}(\{\mathbf{n}_2^{\Lambda} \, \text{odd}\}) \mathbb{P}(\{\mathbf{n}_1^{\g} \equiv 0\}) \mathbb{P}(\{\mathbf{n}_2^{\g} \equiv 0\}).
        \end{split}
    \end{equation}
    Applying Remark \ref{rmk:poisson-rv}, we can rewrite the above expression as
    \begin{equation}
        \begin{split}
            \prod_{x=0}^{N-1} \tanh(J_x) \left[ \mathbb{P}(\{\mathbf{n}^{\Lambda} \, \text{even}\}) \mathbb{P}(\{\mathbf{n}^{\g} \equiv 0\}) \right]^2.
        \end{split}
    \end{equation}

    Additionally, the condition $\partial \mathbf{n} = \varnothing$ is equivalent to $\partial \mathbf{n}^{\Lambda} = \partial \mathbf{n}^{\g}$. By substitution into Equation \eqref{equ:switching-lemma}, we obtain Equation \eqref{equ:cov-identity}.
\end{proof}

To understand the behavior of the denominator of Equation \eqref{equ:cov-identity}, it is worth noting that only $\mathbf{n}^{\g}$ with $\sum_{x = 0}^N \mathbf{n}^{\g}(x)$ being even can contribute to $\mathbb{P}(\{\partial \mathbf{n}^{\Lambda} = \partial \mathbf{n}^{\g}\})$. In the next lemma, we derive a bound on this probability, conditioned on $\sum_{x = 0}^N \mathbf{n}^{\g}(x)$ being even.
\begin{lemma}
    The random current measure on $\{0, \ldots, N\}$ respects the bound
    \begin{equation}
        \frac{\mathbb{P}\left( \{\partial \mathbf{n}^{\Lambda} = \partial \mathbf{n}^{\g}\} \mid \left\{ \sum_{x=0}^N \mathbf{n}^{\g} \, \mathrm{even} \right\} \right)}{\mathbb{P}(\{\mathbf{n}^{\Lambda} \, \mathrm{even}\})} \geq \prod_{x=0}^{N-1} \frac{1+\tanh(J_x)}{2}.
    \end{equation}
    \label{lem:cond-probability-bd-from-below}
\end{lemma}

Because the dependence of $\mathbf{n}^{\Lambda}$ on $\mathbf{n}^{\g}$ under the condition $\{\partial \mathbf{n}^{\Lambda} = \partial \mathbf{n}^{\g}\}$ will be global, we introduce the following terminology:
\begin{definition}
    If $\mathbf{n}^{\g} \in \mathcal{C}^{\g}$, then we will say that an edge $(x,x+1)$ ``splits $\mathbf{n}^{\g}$ evenly'' if $\sum_{y \leq x} \mathbf{n}^{\g}(y)$ and $\sum_{y \geq x+1} \mathbf{n}^{\g}(y)$ are even and that $(x,x+1)$ ``splits $\mathbf{n}^{\g}$ oddly'' if both summations are odd. If $\boldsymbol{\epsilon} = \{\epsilon _x\}_{x=0}^{N-1} \in \{\mathtt{odd}, \mathtt{even}\}^N$, then we will say that $\mathbf{n}^{\g}$ is split according to $\boldsymbol{\epsilon}$ if individual edges $(x,x+1)$ split $\mathbf{n}^{\g}$ evenly whenever $\epsilon _x = \mathtt{even}$ and oddly whenever $\epsilon _x = \mathtt{odd}$.
\end{definition}

The following lemma informs how path-uniqueness simplifies the event $\{\partial \mathbf{n}^{\Lambda} = \partial \mathbf{n}^{\g}\}$.
\begin{lemma}
    Let $\mathbf{n}^{\g} \in \mathcal{C}^{\g}$ be fixed. Then a current $\mathbf{n}^{\Lambda} \in \mathcal{C}^{\Lambda}$ satisfies the relation $\partial \mathbf{n}^{\Lambda} = \partial \mathbf{n}^{\g}$ if and only if $\mathbf{n}^{\g}$ is split evenly by the edges on which $\mathbf{n}^{\Lambda}$ is even and split oddly by the edges on which $\mathbf{n}^{\Lambda}$ is odd. \label{lem:even-odd-splitting}
\end{lemma}

\begin{proof}
    Suppose first that $\partial \mathbf{n}^{\Lambda} = \partial \mathbf{n}^{\g}$. Starting from the left, $\mathbf{n}^{\Lambda}(0)$ must have the same parity as $\mathbf{n}^{\g}(0)$ because no other edges connect to vertex $0$. Proceeding inductively, take $x$ to be such that $\mathbf{n}^{\Lambda}(x)$ has the same parity as $\sum_{y=0}^x \mathbf{n}^{\g}(y)$. If $\mathbf{n}^{\g}(x+1)$ is even, then $\sum_{y=0}^{x+1} \mathbf{n}^{\g}(y)$ will have the same parity as $\sum_{y=0}^x \mathbf{n}^{\g}(y)$, while $\mathbf{n}^{\Lambda}(x+1)$ must have the same parity as $\mathbf{n}^{\Lambda}(x)$ to ensure that an even number of arrivals from $\mathbf{n}^{\Lambda}$ touch $x+1$. On the other hand, if $\mathbf{n}^{\g}(x+1)$ is odd, then $\sum_{y=0}^{x+1} \mathbf{n}^{\g}(y)$ will have the opposite parity as $\sum_{y=0}^j \mathbf{n}^{\g}(y)$, just as $\mathbf{n}^{\Lambda}(x+1)$ and $\mathbf{n}^{\Lambda}(x)$ must have opposite parities so that an odd number of arrivals from $\mathbf{n}^{\Lambda}$ will touch $x+1$. We can thus conclude inductively that $\mathbf{n}^{\Lambda}(x)$ has the same parity as $\sum_{y=0}^x \mathbf{n}^{\g}(y)$ for all $(x,x+1)$. Arguing symmetrically, the same must be true of $\sum_{y=x+1}^N \mathbf{n}^{\g}(y)$.

    Conversely, if the parities of $\mathbf{n}^{\Lambda}(x)$, $\sum_{y=0}^x \mathbf{n}^{\g}(y)$, and $\sum_{y=x+1}^N \mathbf{n}^{\g}(y)$ are the same for all $x \in \{0, \ldots, N-1\}$, we can take any $x \in \Lambda$ and write
    \begin{equation}
        \mathbf{n}^{\g}(x) = \sum_{y = 0}^{N} \mathbf{n}^{\g}(y) - \sum_{y=0}^{x-1} \mathbf{n}^{\g}(y) - \sum_{y=x+1}^{N} \mathbf{n}^{\g}(y),
    \end{equation}
    where on the right-hand side, the first summation is even, the second summation has the same parity as $\mathbf{n}^{\Lambda}(x-1)$ if $x>0$ or is zero if $x=0$, and the third summation has the same parity as $\mathbf{n}^{\Lambda}(x)$ if $x<N$ or is zero of $x=N$. Thus, $\mathbf{n}^{\g}(x)$ has the same parity as the total number of arrivals of $\mathbf{n}^{\Lambda}$ on edges either side of vertex $x$, implying that $x \in \partial \mathbf{n}^{\g}$ if and only if $x \in \partial \mathbf{n}^{\Lambda}$. It follows that $\partial \mathbf{n}^{\Lambda} = \partial \mathbf{n}^{\g}$.
\end{proof}

This characterization lets us prove Lemma~\ref{lem:cond-probability-bd-from-below}.

\begin{proof}[Proof of Lemma~\ref{lem:cond-probability-bd-from-below}]
    Writing 
    \begin{equation}
        \mathbb{P}\left( \{\partial \mathbf{n}^{\Lambda} = \partial \mathbf{n}^{\g}\} \mid \left\{ \sum_{x=0}^N \mathbf{n}^{\g} \, \text{even} \right\} \right) = \mathbb{E}\left[ \mathbb{P}(\{\partial \mathbf{n}^{\Lambda} = \partial \mathbf{n}^{\g}\} \mid \mathbf{n}^{\g})\mid \left\{ \sum_{x=0}^N \mathbf{n}^{\g} \, \text{even} \right\} \right],
    \end{equation}
    Lemma \ref{lem:even-odd-splitting} makes it possible to rewrite the inner conditional probability as
    \begin{equation}
        \begin{split}
            \mathbb{P} (\partial \mathbf{n}^{\Lambda} = \partial \mathbf{n}^{\g} \mid \mathbf{n}^{\g}) = \prod_{x=0}^{N-1} &( \mathbb{P} (\text{$\mathbf{n}^{\Lambda}(x)$ even}) \mathbbm{1}_{\{\text{$(x,x+1)$ splits $\mathbf{n}^{\g}$ evenly}\}} \\&\qquad+ \mathbb{P} (\text{$\mathbf{n}^{\Lambda}(x)$ odd}) \mathbbm{1}_{\{\text{$(x,x+1)$ splits $\mathbf{n}^{\g}$ oddly}\}} ).
        \end{split}
        \label{equ:even-odd-product}
    \end{equation}

    Taking the result of Equation \ref{equ:even-odd-product}, and dividing through by the probability that $\mathbf{n}^{\Lambda}$ is even on all edges, we find
    \begin{equation}
        \frac{\mathbb{P}\left( \{\partial \mathbf{n}^{\Lambda} = \partial \mathbf{n}^{\g}\} \mid \left\{ \sum_{x=0}^N \mathbf{n}^{\g} \, \mathrm{even} \right\} \right)}{\mathbb{P}(\{\mathbf{n}^{\Lambda} \, \mathrm{even}\})} = \mathbb{E} \left[ \prod_{x=0}^{N-1} f_x(\mathbf{n}^{\g}) \mid \left\{ \sum_{x=0}^N \mathbf{n}^{\g}(x) \, \text{even} \right\} \right],
    \end{equation}
    where we have defined
    \begin{align}
        f_x(\mathbf{n}^{\g}) :&= \frac{\mathbb{P}(\{\partial \mathbf{n}^{\Lambda} = \partial \mathbf{n}^{\g}\} \mid \mathbf{n}^{\g})}{\mathbb{P}(\{\mathbf{n}^{\Lambda}(x) \, \text{even}\})} \\
        &= \mathbbm{1}_{\{\text{$(x,x+1)$ splits $\mathbf{n}^{\g}$ evenly}\}} + \tanh(J_x) \mathbbm{1}_{\{\text{$(x,x+1)$ splits $\mathbf{n}^{\g}$ oddly}\}}.
    \end{align}
    
    Due to the result of Lemma \ref{lem:even-odd-splitting}, any $\mathbf{n}^{\g}$ with $\sum_{x=0}^{N} \mathbf{n}^{\g}(x)$ even specifies a unique set of parities $\boldsymbol{\epsilon} \in \{\mathtt{even}, \mathtt{odd}\}^N$ such that $\mathbf{n}^{\g}$ is split according to $\boldsymbol{\epsilon}$, and it can be seen that the definition of $f_x$ for any $x$ only depends on this specified $\boldsymbol{\epsilon}$. To evaluate the conditional expectation, we can consider these functions to take $\{\mathtt{even},\mathtt{odd}\}^N$ as their domain. If we let $\mathtt{odd} \prec \mathtt{even}$ and allow this ordering to extend naturally to a partial ordering on $\{\mathtt{even}, \mathtt{odd}\}^N$, then the set $\{\mathtt{even}, \mathtt{odd}\}^N$ form a distributive lattice on which the functions $f_x$ are a set of non-negative increasing functions. As a result, the FKG inequality guarantees that
    \begin{equation}
        \mathbb{E} \left[ \prod_{x=0}^{N-1} f_{x}(\mathbf{n}^{\g}) \mid \left\{ \sum_{x=0}^{N} \mathbf{n}^{\g}(x) \, \text{even} \right\} \right] \geq \prod_{x=0}^{N-1} \mathbb{E} \left[ f_{x}(\mathbf{n}^{\g}) \mid \left\{ \sum_{x=0}^{N} \mathbf{n}^{\g}(x) \, \text{even} \right\} \right].
    \end{equation}
    Moreover, as Poisson random variables are always more likely to be even than odd, we have that
    \begin{equation}
        \mathbb{P}((x,x+1) \, \text{splits $\mathbf{n}^{\g}$ evenly}) \geq \mathbb{P}((x,x+1) \, \text{splits $\mathbf{n}^{\g}$ oddly}),
    \end{equation}
    implying that
    \begin{equation}
        \mathbb{E} \left[ f_{x}(\mathbf{n}^{\g}) \mid \left\{ \sum_{x=0}^{N} \mathbf{n}^{\g}(x) \, \text{even} \right\} \right] \geq \frac{1 + \tanh(J_x)}{2}.
    \end{equation}
\end{proof}

Combining this bound with the identity given in Lemma~\ref{lem:cov-identity}, we can now prove Theorem~\ref{thm:cov-upper-bd-fm}.
\begin{proof}[Proof of Theorem~\ref{thm:cov-upper-bd-fm}]
    Working from Equation~\eqref{equ:cov-identity}, it suffices to show that
    \begin{equation}
        \frac{\mathbb{P}(\{\mathbf{n}^{\Lambda} \, \mathrm{even}\}) \mathbb{P}(\{\mathbf{n}^{\g} \equiv 0\})}{\mathbb{P}(\{\partial \mathbf{n}^{\Lambda} = \partial \mathbf{n}^{\g}\})} \leq \prod_{x=0}^{N-1} \left( \frac{2}{1+\tanh(J_x)} \right) \frac{1}{\cosh \left( \sum_{x=0}^N h_x \right)}.
        \label{equ:cov-identity-reduced}
    \end{equation}
    To obtain this bound, we can decompose the denominator in the form
    \begin{equation}
        \mathbb{P}(\{\partial \mathbf{n}^{\Lambda} = \partial \mathbf{n}^{\g}\}) = \mathbb{P} \left( \{\partial \mathbf{n}^{\Lambda} = \partial \mathbf{n}^{\g}\} \mid \left\{ \sum_{x=0}^N \mathbf{n}^{\g}(x) \, \text{even} \right\} \right) \mathbb{P} \left( \left\{ \sum_{x=0}^N \mathbf{n}^{\g}(x) \, \text{even} \right\} \right),
    \end{equation}
    recalling that $\{\partial\mathbf{n}^{\Lambda} = \partial \mathbf{n}^{\g}\}$ is only possible if an even number of arrivals of $\mathbf{n}$ touch the ghost site. The bound
    \begin{equation}
        \frac{\mathbb{P}(\{\mathbf{n}^{\Lambda} \, \text{even}\})}{\mathbb{P}\left( \{\partial \mathbf{n}^{\Lambda} = \partial \mathbf{n}^{\g}\} \mid \sum_{x=0}^N \mathbf{n}^{\g} \, \text{even} \right)} \leq \prod_{x=0}^{N-1} \left( \frac{2}{1 + \tanh(J_x)} \right)
    \end{equation}
    is implied by Lemma~\ref{lem:cond-probability-bd-from-below}.  The random variable $\sum_{x=0}^N \mathbf{n}^{\g}(x)$, as the sum of independent Poisson random variables, is itself Poisson distributed with rate $\sum_{x=0}^N h_x$. As noted in Remark~\ref{rmk:poisson-rv},
    \begin{equation}
        \frac{\mathbb{P}(\{\mathbf{n}^{\g} \equiv 0\})}{\mathbb{P}\left( \sum_{x=0}^N \mathbf{n}^{\g}(x) \, \text{even} \right)} = \frac{1}{\cosh \left( \sum_{x=0}^N h_x \right)}.
    \end{equation}
    Multiplying these results yields the inequality \eqref{equ:cov-identity-reduced}.
\end{proof}

\subsection{Proof of Lemma~\ref{lem:cov-signed-and-absolute}}
To prove Lemma~\ref{lem:cov-signed-and-absolute}, we work from results stated in Subsection \ref{ssc:random-currents}.

\begin{proof}
    We recall that the random current measure for coupling constants and external fields $\mathbf{J}$ and $\mathbf{h}$ is identical to the one for parameters $|\mathbf{J}|$ and $|\mathbf{h}|$.  Then, for any $\mathbb{P}$-measurable set $M$ we have
    \begin{equation}
        \left| \mathbb{E} \left[ (-1)^{\#_{\mathbf{J},\mathbf{h}}^-(\mathbf{n}_1 + \mathbf{n}_2)} \mathbbm{1}_{M} \right] \right| \leq \mathbb{E} \left[ \left| (-1)^{\#_{\mathbf{J},\mathbf{h}}^-(\mathbf{n}_1 + \mathbf{n}_2)} \right| \mathbbm{1}_{M} \right] =\mathbb{E} \left[ (-1)^{\#_{|\mathbf{J}|,|\mathbf{h}|}^-(\mathbf{n}_1 + \mathbf{n}_2)} \mathbbm{1}_{M} \right],
    \end{equation}
    as $\#_{|\mathbf{J}|,|\mathbf{h}|}^-(\mathbf{n}_1 + \mathbf{n}_2)$ equals zero. If $M := \{\partial \mathbf{n}_1 = \varnothing; \partial \mathbf{n}_2 = \{i,j\}, \notconnected{i}{\g}{\mathbf{n}_1+\mathbf{n}_2}\}$, Equation \eqref{equ:switching-lemma-alt} can be used to re-express both expressions as
    \begin{equation}
        (\Tilde{Z}_{\mathbf{J},\mathbf{h}} )^2 |\cov_{\mathbf{J},\mathbf{h}} (\sigma_i, \sigma_j)| \leq (\Tilde{Z}_{\mathbf{|J|}, |\mathbf{h}|} )^2 \cov_{\mathbf{|J|},|\mathbf{h}|} (\sigma_i, \sigma_j).
    \end{equation}
    But $\Tilde{Z}_{\mathbf{J},\mathbf{h}} $ and $\Tilde{Z}_{\mathbf{J},|\mathbf{h}|} $ are merely rescalings of $Z_{\mathbf{J},\mathbf{h}} $ and $Z_{\mathbf{J},|\mathbf{h}|} $, respectively, by an identical factor, implying the lemma.
\end{proof}

\subsection{Proof of Theorem \ref{thm:cov-upper-bd-signed}}
Theorem \ref{thm:cov-upper-bd-signed} can now be deduced from Theorem~\ref{thm:cov-upper-bd-fm} and Lemma~\ref{lem:cov-signed-and-absolute}.
\begin{proof}
    As Theorem~\ref{thm:cov-upper-bd-fm} is applicable to the one-dimensional lattice with coupling constants $\mathbf{J}$ and external fields $|\mathbf{h}|$, combining the resulting bound with Lemma~\ref{lem:cov-signed-and-absolute} yields the estimate
    \begin{equation}
        \cov_{\mathbf{J},\mathbf{h}}(\sigma_0, \sigma_N) \leq \left( \prod_{x=0}^{N-1} \frac{4 \tanh(J_x)}{(1+\tanh(J_x))^2} \right) \frac{1}{\cosh^2 \left( {\sum_{x=0}^{N} |h_x|} \right)} \left( \frac{Z_{\mathbf{|J|}, |\mathbf{h}|} }{Z_{\mathbf{J},\mathbf{h}} } \right)^2.
        \label{equ:cov-estimate-partition-ratio}
    \end{equation}
    
    The external field $\mathbf{h}$ can be split according to sign with $h^+_x := \max\{h_x,0\}$ and $h^-_x := |\min\{h_x,0\}|$,
    so that $\mathbf{h} = \mathbf{h}^+ - \mathbf{h}^-$ and $|\mathbf{h}| = \mathbf{h}^+ + \mathbf{h}^-$. Without loss of generality, suppose that $\sum_{x =0}^N h^+_x \geq \sum_{x = 0}^N h^-_x$. We will first show that
    \begin{equation}
        \frac{Z_{\mathbf{J},\mathbf{h}}}{Z_{|\mathbf{J}|,|\mathbf{h}|}} \geq \prod_{x = 0}^N \exp(-2 h^-_x). \label{equ:partition-fn-signed-absolute}
    \end{equation}
    This claim follows from the FKG inequality and the convexity of the exponential function. Observe that
    \begin{equation}
        \frac{Z_{\mathbf{J},\mathbf{h}}}{Z_{|\mathbf{J}|,|\mathbf{h}|}} = \left\langle \prod_{x \in V} \exp({- 2 h_x^- \sigma_x}) \right\rangle_{|\mathbf{J}|,|\mathbf{h}|}.
    \end{equation}
    Because $\exp({-2h^-_x \sigma_x})$ for each $x$ is a positive-valued nonincreasing function of $\boldsymbol{\sigma}$, the FKG inequality is applicable and yields 
    \begin{equation}
        \left\langle \prod_{x \in V} \exp({-2h^-_x \sigma_x}) \right\rangle_{|\mathbf{J}|,|\mathbf{h}|} \geq \prod_{x \in V} \left\langle \exp({-2h^-_x \sigma_x}) \right\rangle_{|\mathbf{J}|,|\mathbf{h}|}.
    \end{equation}
    Jensen's inequality then yields that
    \begin{equation}
        \left\langle \exp({-2h^-_x \sigma_x}) \right\rangle_{\mathbf{J},|\mathbf{h}|} \geq \exp({-2h^-_x \left\langle \sigma_x \right\rangle_{|\mathbf{J}|,|\mathbf{h}|}})
    \end{equation}
    for any $x$, and the fact that $\sigma_x$ cannot exceed $1$ implies the estimate \eqref{equ:partition-fn-signed-absolute}.

    Therefore,
    \begin{align}
        \frac{1}{\cosh^2 \left( {\sum_{x=0}^{N} |h_x|} \right)} \left( \frac{Z_{\mathbf{|J|}, |\mathbf{h}|} }{Z_{\mathbf{J},\mathbf{h}} } \right)^2 &\leq \frac{\prod_{x = 0}^N \exp(-2 h^-_x)}{\cosh \left( {\sum_{x=0}^{N} |h_x|} \right)} \\ 
        &= \frac{2 \exp \left( 2 \sum_{x=0}^N h^-_x \right)}{\exp \left( \sum_{x=0}^N (h^+_x + h^-_x) \right) + \exp \left( - \sum_{x=0}^N (h^+_x + h^-_x) \right)} \\
        &= \frac{2 \exp \left( - \sum_{x=0}^N (h^+_x - h^-_x) \right)}{1 + \exp \left( - 2 \sum_{x=0}^N (h^+_x + h^-_x) \right)} \\
        &= \frac{2 \exp \left( - \sum_{x=0}^N h_x \right)}{1 + \exp \left( - 2 \sum_{x=0}^N |h_x| \right)},
    \end{align}
    and applying this bound to the right-hand side of \eqref{equ:cov-estimate-partition-ratio} completes the proof.
\end{proof}

\section*{Acknowledgements}
JK gratefully acknowledges the financial support of the German Academic Exchange Service (DAAD) as well as the Leibniz Association within the Leibniz Junior Research Group on Probabilistic Methods for Dynamic Communication Networks as part of the Leibniz Competition. 

MG has been supported by NSF awards DMS-1954343 and DMS-2348113.

\appendix

\section{Proof of Lemma~\ref{lem:effective-field}}\label{appendix:truncated-model}

As vertices can be removed iteratively from $\{0, \ldots, N\}$ until the desired sub-lattice is reached, it is sufficient to set $i=0$ and $j=N-1$. We can always find values $A, B \in \mathbb{R}$ such that
\begin{equation}
    \log \left( \sum_{\sigma_N} \exp \left(J_{N-1} \sigma_{N-1} \sigma_N + h_N \sigma_N) \right) \right) = A + B \sigma_{N-1}
\end{equation}
for $\sigma_{N-1} = \pm 1$.

Integrating out the value $\sigma_N$ in the original Gibbs measure $\mu$, we find that, dropping subscripts to indicate coupling constants and external fields,
\begin{align}
    \begin{split}
        \mu(\boldsymbol{\sigma}|\{0, \ldots, N-1\}) &= \frac{1}{Z} \sum_{\sigma_N} \exp \left( {\sum_{x=i}^{j-1} (J_x \sigma_x \sigma_{x+1} + h_x \sigma_x) + h_j \sigma_j} \right) \\&\qquad\qquad\;\cdot \, \exp(J_{N-1} \sigma_{N-1} \sigma_N + h_N \sigma_N)
    \end{split} \\
    &= \frac{1}{Z} \exp \left( {\sum_{x=0}^{N-1} (J_x \sigma_x \sigma_{x+1} + h_x \sigma_x) + A + h_{N-1}' \sigma_{N-1}} \right),
\end{align}
where we have defined $h'_{N-1} := h_{N-1} + B$.

In order to bound the correction $B$, we note that
\begin{align}
    B &= \frac{1}{2} \log \left( \frac{\sum_{\sigma_N} \exp \left(J_{N-1} \sigma_N + h_N \sigma_N) \right)}{\sum_{\sigma_N} \exp \left(-J_{N-1} \sigma_N + h_N \sigma_N) \right)} \right) \\
    &= \frac{1}{2} \log \left( \frac{\sum_{\sigma_N} \exp \left(-J_{N-1} \sigma_N + h_N \sigma_N) \right) \exp(2J_{N-1} \sigma_N)}{\sum_{\sigma_N} \exp \left(-J_{N-1} \sigma_N + h_N \sigma_N) \right)} \right) \\
    &\leq \frac{1}{2} \log \left( \frac{\exp(2J_{N-1}) \sum_{\sigma_N} \exp \left(-J_{N-1} \sigma_N + h_N \sigma_N) \right)}{\sum_{\sigma_N} \exp \left(-J_{N-1} \sigma_N + h_N \sigma_N) \right)} \right) \\
    &= J_{N-1},
\end{align}
and that we can similarly establish that $B \geq - J_{N-1}$, so that
\begin{align}
    |h'_{N-1} - h_{N-1}| &\leq J_{N-1}.
\end{align}

\bibliography{bib}
\bibliographystyle{alpha}

\end{document}